\documentclass{article}
\usepackage{amsmath}
 \usepackage{amsfonts}
\begin{document}
\title{{Singular solutions of conformal Hessian  equation}} \author{{Nikolai Nadirashvili\thanks{I2M, Aix-Marseille Universit\'e, 39, rue F. Joliot-Curie, 13453
Marseille  FRANCE, nicolas@cmi.univ-mrs.fr},\hskip .4 cm Serge
Vl\u adu\c t\thanks{I2M, Aix-Marseille Universit\'e, Luminy, case 907, 13288 Marseille Cedex
FRANCE and IITP RAS, B.Karetnyi,9, Moscow, RUSSIA,
vladut@iml.univ-mrs.fr} }}

 \date{}
\maketitle
\def\S{\mathbb{S}}
\def\Z{\mathbb{Z}}
\def\R{\mathbb{R}}
\def\C{\mathbb{C}}
\def\N{\mathbb{N}}
\def\H{\mathbb{H}}
\def\O{\mathbb{O}}
\def\tilde{\widetilde}

\def\n{\hfill\break} \def\al{\alpha}\def\epsilon{\varepsilon} \def\be{\beta} \def\ga{\gamma} \def\Ga{\Gamma}
\def\om{\omega} \def\Om{\Omega} \def\ka{\kappa} \def\lm{\lambda} \def\Lm{\Lambda}
\def\dl{\delta} \def\Dl{\Delta} \def\vph{\varphi} \def\vep{\varepsilon} \def\th{\theta}
\def\Th{\Theta} \def\vth{\vartheta} \def\sg{\sigma} \def\Sg{\Sigma}
\def\bendproof{$\hfill \blacksquare$} \def\wendproof{$\hfill \square$}
\def\holim{\mathop{\rm holim}} \def\span{{\rm span}} \def\mod{{\rm mod}}
\def\rank{{\rm rank}} \def\bsl{{\backslash}}
\def\il{\int\limits} \def\pt{{\partial}} \def\lra{{\longrightarrow}}

{\em Abstract.} 
We show that for any $\epsilon\in ]0,1[$ there exists an  analytic outside zero solution to a uniformly elliptic conformal Hessian equation in a ball $B\subset\R^5$ which belongs to $C^{1,\epsilon} (B)\setminus C^{1,\epsilon+} (B)$.
\bigskip



\section{Introduction}
\bigskip

In this paper  we study a class of fully nonlinear second-order elliptic equations of the form
\begin{equation}F(D^2u, Du,u)=0\end{equation}
defined in a domain of ${ \R}^n$. Here $D^2u$ denotes the
Hessian of the function $u,\;Du$ being its gradient. We assume that
$F$ is a Lipschitz  function defined on a domain in  the space $ {\rm Sym}_2({ \R}^n)\times { \R}^n\times { \R},\;
 {\rm Sym}_2({ \R}^n)$ being the space of ${n\times n}$  symmetric matrices and that  $F$  satisfies the uniform 
 ellipticity condition, i.e. there exists a constant $C=C(F)\ge 1$ (called an {\it ellipticity constant}) such that 
$$C^{-1}||N||\le F(M+N)-F(M) \le C||N||\;$$ 
for any non-negative definite symmetric matrix $N$; if
$F\in C^1({\rm Sym}_2({ \R}^n))$ then this condition is equivalent to
$$\frac{1}{ C'}|\xi|^2\le F_{u_{ij}}\xi_i\xi_j\le C' |\xi |^2\;,\forall\xi\in { \R}^n\;. $$  
Here, $u_{ij}$ denotes the partial derivative $\pt^2 u/\pt x_i\pt x_j$. A function $u$ is called a {\it classical} 
solution of (1) if $u\in C^2(\Om)$ and $u$ satisfies (1).  Actually, any classical solution of (1) is a smooth 
($C^{\alpha +3}$) solution, provided that $F$ is a smooth $(C^\alpha )$ function of its arguments.

 More precisely, we are interested in conformal Hessian  equations (see, e.g. \cite{T},  pp. 5-6) i.e. those of the form 
\begin{equation} F[u] := f(\lambda(A^u)) = \psi( u, x) \end{equation}
$f$   being  a  Lipschitz  function on ${ \R}^n$ invariant under  permutations of the coordinates and
 $$\lambda(A^u)=(\lambda_1, \ldots, \lambda_n)$$ 
being the eigenvalues of the conformal Hessian in ${\R}^n$:
\begin{equation} \label{confhess}A^u:=uD^2u-\frac{1}{2}|Du|^2I_n \end{equation}
where $n\ge 3, u>0$.

In this case $F$ is invariant under conformal mappings $T :{\R}^n \longrightarrow{\R}^n,$ i.e. transformations which preserve angles between curves. In contrast to the case $n=2$, for $n\ge 3$ any conformal transformation  of ${\R}^n $ 
is decomposed into a finitely many  M\"obius transformations, that is mappings of the form
$$ Tx = y +\frac{kA(x - z)} {  |x-z|^a}, $$
with $x,z\in {\R}^n, k\in {\R}, a\in \{0,2\}$ and an orthogonal matrix $A$. In other words, each $T$ is a composition 
of a translation, a homothety,  a rotation and (may be)  an inversion. If $T$ is a conformal mapping and $v(x) = J_T^{-1/n} u  (T x)$, where $J_T$ denotes the Jacobian determinant of $T$ then $F[v] = F[u].$ Note that this class of equations is very important in geometry, see \cite{li} and references therein.

\medskip
We are interested in the Dirichlet problem
\begin{equation}\label{dir}\begin{cases}F(D^2u, Du,u)=0, u>0 &\text{in}\;
 \Om \cr \quad \quad  u=\vph & \text{on}\; \pt\Om\;,\cr\end{cases}\end{equation}
 where  $\Omega \subset {\R}^n$ is a bounded domain with a smooth boundary $\partial \Omega$
and $\vph$ is a continuous function on $\pt\Om$.

Consider the problem of existence and regularity of solutions to the Dirichlet problem (4)  which has always a unique viscosity (weak) solution for  fully nonlinear elliptic equations. The viscosity solutions  satisfy the equation
(1) in a weak sense, and the best known interior regularity (\cite{C},\cite{CC},\cite{T3}) for them is $C^{1+\epsilon }$ for some $\epsilon> 0$. For more details see \cite{CC}, \cite{CIL}. Recall that in  \cite{NV1} the authors constructed a homogeneous  singular  viscosity solution in 5 dimensions for Hessian equations of  order $1+\delta$ for any  $\delta  \in ]0,1]$, that is, of any order compatible with  the mentioned interior regularity results.
In fact we proved in \cite{NV1} the following result.\pagebreak
 
\medskip{\bf Theorem 1.1.}

{\it The function
$$w_{5,\delta}(x)=P_{5} (x)/ |x|^{1+\delta }, \;\delta\in [0,1[
$$ is a viscosity solution  to  a uniformly elliptic Hessian equation
  $F(D^2w)=0$ with a smooth functional $F$ in a unit ball $B\subset {\R}^{5}$  for  the isoparametric Cartan cubic form
$$ P_{5}(x)=x_1^3+\frac{3 x_1}2\left(z_1^2 + z_2^2-2 z_3^2-2x_2^2\right)+\frac{3\sqrt 3}2\left(x_2z_1^2-x_2z_2^2 + 2z_1z_2z_3\right)$$
with  $x=(x_1,x_2,z_1,z_2,z_3)$.}\medskip

\noindent which proves  the optimality of  the interior $C^{1+\epsilon}$-regularity of viscosity solutions to fully nonlinear equations  in 5  and more dimensions.

In the present paper we show that the same singularity result remains true for conformal Hessian equations.

\medskip{\bf Theorem 1.2.}

{\it Let  $\delta\in ]0,1[.$ The function 
$$u(x):=c+w_{5,\delta}(x)=c+\frac{P_{5} (x)}{ |x|^{1+\delta }}, $$ 
is a viscosity solution  to  a uniformly elliptic conformal Hessian equation   $(1)$ in a unit ball $B\subset {\R}^{5}$  
 for a sufficiently large positive constant $c\; (c=240000 $ is sufficient for $\delta=\frac 12).$}\medskip

Notice also that the result does {\em not} hold for $\delta=0$ and we do not know how to construct a non-classical $C^{1,1}$-solution   
to  a uniformly elliptic conformal Hessian equation.

 The rest of the paper is organized as follows: in Section 2 we recall some necessary preliminary results   and we prove our main results in Section 3; to simplify  the notation  we suppose that $\delta=\frac 12$ in Section 3; for any $\delta$ the  proof is along the same line, but more cumbersome.
The proof in Section 3   uses MAPLE to varify some algebraic identities but is completely rigorous (and is human-controlled for  $\delta=\frac 12$).

\section{Preliminary results }

{\em Notation}: for a real  symmetric matrix $A$ we denote by  $|A|$ the maximum of the absolute value of its eigenvalues.

Let $ u$ be a strictly positive function on $B_1 .$ Define the map
$$\Lambda :B_1  \longrightarrow \lambda (S) \in { \R}^n\; .$$
$\lambda(S)=\{  \lambda_1\ge...\ge\lambda_n\} \in { \R}^n$ being  the (ordered) set  of eigenvalues of the conformal Hessian 
$$A^u:=uD^2u-\frac{1}{2}|Du|^2I_n.$$

The following ellipticity criterion can proved similarly to  Lemma 2.1 of \cite{NV4}.

\bigskip {\bf  Lemma 2.1.} {\it Suppose that the family 
$$\{A^u(a)-O^{-1}\cdot A^u(b)\cdot  O\: : a,b\in B_1, O\in {\mathrm SO}(n)\}\setminus\{0\}$$ is uniformly hyperbolic, i.e. if $\{\mu_1(a,b,O)\ge \ldots \ge \mu_n(a,b,O)\}$ is the ordered spectrum of  $A^u(a)-{O^{-1}}\cdot A^u(b)\cdot O\neq 0$ then 
$$\forall  a,b\in B_1,\forall O\in {\mathrm SO}(n), \:\;\;C^{-1}\le -\frac{\mu_1(a,b,O)}{\mu_n(a,b,O)}\le C $$  for some constant $C>1$.
   Then $u$ is a viscosity solution in $B_1$ of a uniformly elliptic conformal Hessian equation $(1)$.}\medskip

We recall then some properties of the function $w:=w_{5,\delta}(x)=\frac{P_{5} (x)}{ |x|^{1+\delta }}, $   
and its Hessian  $D^2w$ proved in  \cite{NV1}.\medskip
 
\medskip{\bf Lemma 2.2.}

 {\em There exists a 3-dimensional  Lie subgroup $G_P$ of ${\rm SO}(5)$ such that $P$ is invarant under its natural action 
and  the orbit $G_P\S_1^1$ of the circle
 $$\S_1^1=\{(\cos(\chi),0,\sin(\chi),0,0):\chi\in { \R}\}\subset \S_1^4  $$
under this action is the whole $\S_1^4.$}
 \medskip

  \medskip{\bf Lemma 2.3. }

\noindent $(i)$\hskip .5 cm {Let \em $x\in \S_1^4$, and let   $x\in G_P(p,0,r,0,0)$ with $p^2+r^2=1$.
Then $$Spec(D^2w_{5,\delta}(x))= \{\mu_{1,\delta},\mu_{2,\delta},\mu_{3,\delta},\mu_{4,\delta},\mu_{5,\delta} \}$$
for $$\mu_{1,\delta}= \frac{p(p^2\delta+6-3\delta)}{2},$$   
$$\mu_{2,\delta}= \frac{p(p^2\delta-3-3\delta)+3\sqrt{12-3p^2}}{2}, $$
$$\mu_{3,\delta}=\frac{p(p^2\delta-3-3\delta)-3\sqrt{12-3p^2}}{2}, $$ 
$$\mu_{4,\delta}= -\frac{p\delta(6-\delta)(3-p^2)+\sqrt{D(p,\delta)}}{4},$$
$$\mu_{5,\delta}= -\frac{p\delta(6-\delta)(3-p^2)-\sqrt{D(p,\delta)}}{4} ,$$
and
$$D(p,\delta):=(6-\delta)(4-\delta)(2-\delta)\delta(p^2-3)^2p^2+144(\delta-2)^2>0. $$}

\noindent $(ii)$\hskip .5 cm{\em Let  $\lambda_1\ge\lambda_2\ge\ldots\ge\lambda_5$ be the ordered eigenvalues  of $D^2w_{5,\delta}(x)$.
Then
$$\lambda_1= \mu_{2,\delta}, \quad \lambda_5= \mu_{3,\delta},\quad\quad\quad\quad$$
\begin{equation*}\lambda_2 = \begin{cases} \mu_{4,\delta} & \text{for}\; p\in [-1,p_0(\delta)],\cr
\mu_{1,\delta} & \text{for}\; p\in [p_0(\delta),1],  \end{cases}\end{equation*}

\begin{equation*}\lambda_3 = \begin{cases} \mu_{5,\delta} & \text{for}\; p\in [-1,-p_0(\delta)],\cr
\mu_{1,\delta} & \text{for}\; p\in [-p_0(\delta),p_0(\delta)], \cr 
\mu_{4,\delta} & \text{for}\; p\in [ p_0(\delta),1], \end{cases}\end{equation*}
\begin{equation*}\lambda_4 = \begin{cases} \mu_{1,\delta} & \text{for}\; p\in [-1,-p_0(\delta)],\cr\mu_{5,\delta} & \text{for}\; p\in [-p_0(\delta),1],\cr \end{cases}\end{equation*}
 where $$p_0(\delta):=\frac{3^{1/4}\sqrt{1-\delta}}{(3+2\delta-\delta^2)^{1/4}}=\frac{3^{1/4}\sqrt{ \varepsilon}}{(4 - \varepsilon^2)^{1/4}}\in ]0,1].$$}
\medskip

Note  the oddness property of the spectrum:
$$\lambda_{1,\delta}(-p)=-\lambda_{5,\delta}(p),\;\lambda_{2,\delta}(-p)=-\lambda_{4,\delta}(p),
\;\lambda_{3,\delta}(-p)=-\lambda_{3,\delta}(p).$$

\medskip{\bf Proposition 2.1.} 

{\em Let  $N_\delta(x)={D^2w_{\delta}}(x), \;   0\leq\delta <1$. Suppose that  $a\neq b \in B_1\setminus\{0\}$  and let 
$O\in {\hbox {O}}({5} )$ be an orthogonal matrix s.t.
 $$N_\delta(a,b,O):=N_\delta(a)- {^tO}\cdot N_\delta(b)\cdot O\neq 0.$$
 Denote $ \Lambda_1\ge\Lambda_2\ge \ldots\ge\Lambda_{5}$  the eigenvalues of the matrix
 $N_\delta(a,b,O).$
  Then
$$   \frac{1}{C}\le - \frac{\Lambda_1}{\Lambda_{5}}\le C   $$
for $C:=C(\delta):=  \frac{1000(\delta+1)(3-\delta)}{3(1-\delta)^2};$ for $k\in [\frac{1}{2},1]$ one can choose $C=1000$. } 

\medskip{\bf Corollary 2.1.} $$ \Lambda_1\ge \frac{|N_\delta(a,b,O)|}{C(\delta)},\:  |\Lambda_5|\ge \frac{|N_\delta(a,b,O)|}{C(\delta)}.$$

\medskip We need also the following classical Weyl's  result:

 \medskip {\bf  Lemma 2.4.} 
  
{\em Let $ A,B$  be two real symmetric  matrices with the eigenvalues $\lambda_1\ge\lambda_2\ge\ldots\ge\lambda_{n} $ and 
$\lambda'_1\ge\lambda'_2\ge\ldots\ge\lambda'_{n}$ respectively.
Then for the eigenvalues $\Lambda_1\ge\Lambda_2\ge\ldots\ge\Lambda_{n} $ of the matrix $A-B$ we have}
$$ \Lambda_1\ge\max_{i=1,\cdots, n}(\lambda_i-\lambda'_i), \;\;\Lambda_n\le\min_{i=1,\cdots, n}(\lambda_i-\lambda'_i).$$

\section{Proofs} Let $n=5, u(x)=c+w_{5,\delta}(x)$. We begin with $\delta=0 $ and show that the result is false in this case. Indeed let $a=(1,0,0,0,0), b=(\frac 12,0,0,0,0), O=I_5.$ Then
 $$w(a)=1, w(b)=\frac 12, |Du(a)|=|Dw(a)|=9,  |Du(b)|^2=|Dw(b)|^2=\frac 94,$$
 $$ D^2u(a)= D^2w(a)= D^2u(b)= D^2w(b),$$ 
and $$A^u(a)-  A^u(b) =\frac 12 D^2w(a)-\frac{27}{4}I_5$$
  which is negative since the spectrum of $D^2w(a)$ is $( 2, 2, 2,-7, -7).$ The reason is clearly that  $D^2w(a)$ for $\delta=0 $ is homogeneous order 0 and   depends only  on the direction vector $a/|a|$.

Suppose now that $\delta\in ]0,1[$. As we mentioned before, we set $\delta=\frac 12;$ in this case $c=240000$. First we spell out Lemma  2.3 for  $\delta=\frac 12.$

 \medskip {\bf Lemma 3.1. }

\noindent $(i)$\hskip .5 cm {Let \em $x\in \S_1^4$, and let   $x\in G_P(p,0,r,0,0)$ with $p^2+r^2=1$.
Then $$Spec(D^2u (x))= Spec(D^2w (x))= \{\mu_{1 },\mu_{2 },\mu_{3 },\mu_{4 },\mu_{5 } \}$$
for $$\mu_{1 }= \frac{3p(p^2+1)}{4},  $$   
$$\mu_{2 }= \frac{ 3p(p^2-5)+6\sqrt{12-3p^2}}{4}, $$
$$\mu_{3 }=\frac{ 3p(p^2-5)-6\sqrt{12-3p^2}}{4}, $$ 
$$\mu_{4 }= \frac {27p(p^2-3)+3\sqrt{105p^6-630p^4+945p^2+64}}{16},$$
$$\mu_{5 }=  \frac {27p(p^2-3)-3\sqrt{105p^6-630p^4+945p^2+64}}{16} .$$
  
\noindent $(ii)$\hskip .5 cm  Let  $\lambda_1\ge\lambda_2\ge\ldots\ge\lambda_5$ be the ordered eigenvalues  of $Spec(D^2u (x))= Spec(D^2w (x))$.
Then
$$\lambda_1= \mu_{2}, \quad \lambda_5= \mu_{3},\quad\quad\quad\quad$$
\begin{equation*}\lambda_2 = \begin{cases} \mu_{4} & \text{for}\; p\in [-1,p_0],\cr
\mu_{1} & \text{for}\; p\in [p_0,1],  \end{cases}\end{equation*}

\begin{equation*}\lambda_3 = \begin{cases} \mu_{5} & \text{for}\; p\in [-1,-p_0],\cr
\mu_{1} & \text{for}\; p\in [-p_0,p_0], \cr 
\mu_{4} & \text{for}\; p\in [ p_0,1], \end{cases}\end{equation*}
\begin{equation*}\lambda_4 = \begin{cases} \mu_{1} & \text{for}\; p\in [-1,-p_0],\cr\mu_{5} & \text{for}\; p\in [-p_0,1],\cr \end{cases}\end{equation*}
 where}$$p_0 = 5^{-1/4}\simeq 0.6687403050.$$
\medskip
We will need also the derivatives of the egenvalues.

 \medskip
{\bf Lemma 3.2. }  {Let \em $d_i(p):=\frac{d(\mu_i)}{dp}$. Then }
 $$d_{1 }(p)= \frac{3(3p^2+1)}{4},  $$
$$d_{2 }(p)= -\frac{ 3(5-3p^2)}{4}+ \frac{9p}{2\sqrt{12-3p^2}}, $$
$$d_{3 }(p)= -\frac{ 3(5-3p^2)}{4}- \frac{9p}{2\sqrt{12-3p^2}}, $$
 $$d_{4 }(p)=\frac{81 (1-p^2)}{16}\left(\frac{35p(3-p^2)}{3\sqrt{105p^6-630p^4+945p^2+64}} -1 \right),$$
 $$d_{5 }(p)=-\frac{81 (1-p^2)}{16}\left(\frac{35p(3-p^2)}{3\sqrt{105p^6-630p^4+945p^2+64}}  +1\right).$$

Simple calculus gives

 \medskip{\bf Corollary 3.1. } $$D:= \max\{ \left|d_i(p)\right|:p\in [-1,1],\; i=1,\ldots,5\}<10.$$ 

Below we denote  $D_i(p):=\frac{d(\lambda_i)}{dp}$; the relation of $D_i(p)$ and $d_i(p)$ is clear from Lemma 3.1 (ii);
for example, $D_1(p)=d_2(p),\:D_5(p)=d_3(p).$ \medskip

The proof of Theorem 1.2 is based on the following lemmas. Let $$  a,b\in B_1\setminus\{0\}, |a|=s\le 1, |b|=t\le 1, O\in {\hbox {O}}({5} ),$$
$$a':=\frac as \in  G_P(p,0,r,0,0),\: b':=\frac bt\in G_P(q,0,r',0,0).$$ Below we denote 
$$K:=K(p,q,s,t)= |s-t|+| p-q|,$$  
$$M_1:=M_1(a,b,O):=D^2u(a) -{O^{-1}}D^2u(b)\cdot O,$$
$$ M_2:=M_2(a,b,O):=w(a)D^2u(a) -{O^{-1}}w(b)D^2u(b)\cdot O.$$

\medskip{\bf Lemma 3.3. }  $$ \left|\left|Du(a)\right|^2 -\left|Du(b)\right|^2\right|\le 16K .$$

\smallskip{\em Proof.} First, $|Du(a)|^2=|Dw(a)|^2,  |Du(b)|^2 =|Dw(b)|^2.$ Since $P=P_5(x)$ can be represented as the generic traceless norm 
in the Jordan algebra ${\rm Sym}_3(\R)$ it verifies the eiconal equation $|D P|^2=|x|^4$, see e.g. \cite{Tk}. Therefore, an easy calculation gives
 $$|Du(a)|^2=\frac{9s(16-3p^2(p^2-3)^2)}{32},\; |Du(b)|^2=\frac{9t(16-3q^2(q^2-3)^2)}{32},$$
 $$ \left|\left|Du(a)\right|^2 -\left|Du(b)\right|^2\right|\le \left|\frac{9s(16-3p^2(p^2-3)^2)}{32}-\frac{9t(16-3p^2(p^2-3)^2)}{32}\right|+$$
$$ +\left|\frac{9t(16-3p^2(p^2-3)^2)}{32}-\frac{9t(16-3q^2(q^2-3)^2)}{32}\right|=$$
$$ \left|\frac{9(s-t)(16-3p^2(p^2-3)^2)}{32}\right|+\left|\frac{27t(p-q)(p+q)((q^2-3)^2-(p^2-3)^2)}{32}\right|\le$$
$$ \left|\frac{9(s-t)}{2}\right|+\left|\frac{243(p-q)}{16}\right|\le 16K.$$

\medskip{\bf Lemma 3.4. }  {\em Let $M:=\left|M_1\right|=\left| D^2u(a) -{O^{-1}}\cdot D^2u(b)\cdot O\right|$. Then 
  $$  M\ge\frac K8.$$}
{\em Proof.}  If one replaces $a$ by $a'=a/s$ and $b$ by $b''=b/s$ the quantity $M$ gets bigger and $K$ gets smaller.
Therefore we can suppose that $|a|=s=1$. Then we have
  $$D^2u(a) -{O^{-1}}\cdot D^2u(b)\cdot O=D^2u(a) -\frac{{O^{-1}}\cdot D^2u(b')\cdot O}{\sqrt t}.$$

By  Lemma 2.4 we have 
$$M\ge \max\left\{\lambda_i(p)-  \frac{\lambda_i(q)}{\sqrt t}:i=1,\ldots,5\right\},$$
$$M\ge  \left|\min\left\{\lambda_i(p)-  \frac{\lambda_i(q)}{\sqrt t}:i=1,\ldots,5\right\}\right|.$$
Suppose first $p\ge q$. If $q\ge -\frac{24}{25}=-0.96$ then  
$$\forall p'\in [q, p],\:D_1(p')<-1/4=-0.25,\:\lambda_1(p)>\frac 32 $$
 (by a simple calculation using the explicit formulas for $D_1,\lambda_1)$. Therefore 
$$\lambda_1(p)-  \frac{\lambda_1(q)}{\sqrt t}=\lambda_1(p)-\lambda_1(q)+\lambda_1(q)- \frac{\lambda_1(q)}{\sqrt t}\le -\frac{p-q}{4}-\frac 32 \left( \frac{1}{\sqrt t}-1\right)< -\frac K4.$$
If $q<-0.96$ but $p\ge  -\frac{23}{25}=-0.92$ then 
$$\lambda_1(p)- \frac{\lambda_1(q)}{\sqrt t}=\lambda_1(p)-\lambda_1(q)+\lambda_1(q)-\frac{\lambda_1(q)}{\sqrt t} \le\lambda_1(p)-\lambda_1 \left(\frac{24}{25}\right)+\lambda_1(q)- \frac{\lambda_1(q)}{\sqrt t}$$ 
$$ -\frac{p+0.96}{4}-\frac 32 \left( \frac{1}{\sqrt t}-1\right)< -\frac{p-q}{8}-\frac 32 \left( \frac{1}{\sqrt t}-1\right)< -\frac K8.$$
Suppose then that $q<-0.96,\: p<-0.92.$ In this case  we have 
$$\forall p'\in [q,p] ,\; d_2(p')>\frac 52,\; \lambda_2(p')<-\frac 32$$
 and thus
$$\lambda_2(p)-  \frac{\lambda_2(q)}{\sqrt t}=\lambda_2(p)-\lambda_2(q)+\lambda_2(q)- \frac{\lambda_2(q)}{\sqrt t}\ge \frac{5(p-q)}{2}+\frac 32 \left( \frac{1}{\sqrt t}-1\right)\ge\frac{3 K}{4}$$ which finishes the proof for $p\ge q$. The case $q\ge p$ is treated  similarly
(replace $\lambda_1$ by $\lambda_5$ and $\lambda_2$ by $\lambda_4$).

\medskip{\bf Lemma 3.5. }$$\left| M_2\right| = \left|w(a)D^2u(a) -{O^{-1}}w(b)D^2u(b)\cdot O\right|\le 10K$$

\smallskip{\em Proof.}  Indeed, let $a':=a/s,b':=b/s$ then by homogeneity $$ \left|w(a)D^2w(a) -{O^{-1}}w(b)D^2w(b)\cdot O\right|= \left| s D^2w(a') -{O^{-1}}\cdot t D^2w(b')\cdot O\right|\le $$
$$ \le s\left|  D^2w(a') -{O^{-1}}\cdot  D^2w(b')\cdot O\right|+|s-t| \cdot |{O^{-1}}\cdot D^2w(b')|\le$$
$$\le\max_{p,i}\{ |D_i(p)|\}|p-q|+7|s-t|= \max_{p,i}\{ |d_i(p)|\}|p-q|+7|s-t|\le 10K.$$

\medskip{\em Remark 3.1. } These results remain true for any $\delta\in ]0,1[$ if one replaces the respective  constants 16, 1/8 and 10 in Lemmas 3.3, 3.4 and 3.5 by appropriate positive constants depending on $\delta$. On the contrary, Lemma 3.4 is false for $\delta=0$.\medskip

We can now prove the uniform hyperbolicity of  $M(a,b,O)$ and thus the theorem. In fact we show that one can take $C=6007$ in Lemma 3.1.

Indeed,  $$|M(a,b,O)|=\left|A^u(a)-{O^{-1}}\cdot A^u(b)\cdot O\right| =\left|cM_1+M_2-\left(|Du(a)|^2 -|Du(b)|^2\right)I_5\right|.$$ 
 Therefore, $$|\Lambda_5|\ge c|\Lambda_5(M_1)|-10K-16K\ge \frac{ c|M_1|}{1000}-26K\ge 240|M_1|-26K\ge 4K,$$
 $$|\Lambda_1|\ge c\Lambda_1(M_1)-10K-16K\ge \frac{ c|M_1|}{1000}-26K\ge 240|M_1|-26K\ge 4K,$$
$$|M(a,b,O)|\le c\left|M_1\right|+|M_2|+\left||Du(a)|^2 -|Du(b)|^2\right|\le  c\left|M_1\right|+26K.$$
 Thus 
   $$\frac 1C< \frac{4}{240026}\le\frac{240|M_1|-26K}{c\left|M_1\right|+26K}\le\frac{|\Lambda_5|}{|\Lambda_1|}\le\frac{c\left|M_1\right|+26K}{240|M_1|-26K}\le\frac{240026}{4}< C$$ 
which finishes the proof. Notice that we can take $C=1000+\epsilon$ for $\delta\le \frac 12$ if $c$ is sufficiently large;
in the case $\frac 12<\delta< 1$ for sufficiently large $c$ one gets $C=C(\delta)+\epsilon=\frac{1000(\delta+1)(3-\delta)}{3(1-\delta)^2}+\epsilon.$

 \end{document}